\documentclass[11pt]{amsart}

\usepackage{amsmath,amsfonts,amsthm,amssymb, mathtools}
\usepackage{setspace}
\usepackage{ mathrsfs }
\usepackage{ marvosym }

\usepackage{Tabbing}
\usepackage{fancyhdr}
\usepackage{lastpage}
\usepackage{extramarks}

\usepackage{chngpage}
\usepackage{soul,color}
\usepackage{graphicx,float,wrapfig}
\usepackage{fullpage}
\usepackage{hyperref}
\usepackage[all]{xy}

\usepackage{enumerate}
\usepackage{marvosym}
\usepackage{enumitem}
\newcommand{\be}{\begin{enumerate}}
\newcommand{\ee}{\end{enumerate}}
\newcommand{\beq}{\begin{equation}}
\newcommand{\eeq}{\end{equation}}
\usepackage{verbatim}

\newcommand{\ba}{\begin{align*}}

\newcommand{\ea}{\end{align*}}

\newcommand{\D}{{\mathbb D}}
\newcommand{\R}{{\mathbb R}}

\newcommand{\T}{{\mathcal{T}}}

\newcommand{\C}{{\mathbb C}}

\renewcommand{\H}{{\mathbb H}}

\newcommand{\PSL}{\operatorname{PSL}}

\newcommand{\g}{\mathbf{g}}

\newcommand{\im}{\text{Im}}

\newcommand{\Aut}{\text{Aut}}

\newcommand{\stab}{\text{Stab}}

\newcommand{\id}{\text{id}}
\newcommand{\ol}[1]{\overline{#1}}

\newcommand{\pd}[2]{\frac{\partial#1}{\partial#2}}

\newcommand{\ra}{\rightarrow}
\newcommand{\tra}\twoheadrightarrow

\newcommand{\tla}\twoheadleftarrow

\newcommand{\vn}{\varepsilon}

\newcommand{\lp}{\left(}

\newcommand{\rp}{\right)}
\newcommand{\lpi}{\left|}
\newcommand{\rpi}{\right|}

\newcommand{\lbr}{\left\{}
\newcommand{\rbr}{\right\}}

\newcommand{\lbrac}{\left[}

\newcommand{\rbrac}{\right]}

\newcommand{\z}{\mathbf{z}}

\def\bea#1\eea{\begin{align*}#1\end{align*}}
\def\bc#1\ec{\begin{comment}#1\end{comment}}

\newtheorem{Theorem}{Theorem}[section]

\newtheorem{Lemma}[Theorem]{Lemma}
\newtheorem{Proposition}[Theorem]{Proposition}
\newtheorem{Corollary}[Theorem]{Corollary}

\let\lim\relax \DeclareMathOperator*\lim{lim\vphantom{p}}

\title{Asymptotics of the translation flow on holomorphic maps out of the poly-plane}
\date{}
\author{Dmitri Gekhtman}

\begin{document}

\maketitle

\begin{abstract}
We study the family of holomorphic maps from the polydisk to the disk which restrict to the identity on the diagonal.
In particular, we analyze the asymptotics of the orbit of such a map under the conjugation action of a unipotent subgroup
of $\PSL_2(\R)$. We discuss an application of our results to the study of the Carath\'eodory metric on Teichm\"uller space. 
\end{abstract}

\section{Introduction}
Let $\H$ be the upper half-plane $\H = \{\lambda\in \C| \im(\lambda)>0\}$. The poly-plane $\H^n= \H\times \cdots \times \H$ is the $n$-fold product of $\H$ with itself. 

Let $\mathcal{D}$ be the family of holomorphic functions $f:\H^n \ra \H$ which restrict to the identity on the diagonal,
i.e. $f(\lambda,\ldots, \lambda)=\lambda$ for all $\lambda\in \H$.
Fix $t\in \R$. If $f$ is in $\mathcal{D}$, then so is the map $f_t$ defined by
\begin{equation}\label{tflow}
f_t(z_1,\ldots, z_n) = f(z_1-t,\ldots, z_n-t) + t.
\end{equation}
The action $(f,t)\mapsto f_t$ is called the {\em translation flow} on $\mathcal{D}$.

In this paper, we study the asymptotics of the translation flow.
Suppose $f \in \mathcal{D}$, and let $\alpha_j = \pd{f}{z_j}(i,\ldots, i)$ for $j=1,\ldots, n$.
Our main result is that for ``most" $t\in \R$,  $f_t$ is ``close" to the translation-invariant function $\g(z_1,\ldots, z_n) = \sum_{j=1}^n \alpha_j z_j.$
More precisely, we prove
%  \vspace{2mm}\\
%{\bf Theorem \ref{Main}}. {\em 
%Fix $\vn>0$, and let $U$ be any open neighborhood of $\g$ in the compact-open topology. Then for sufficiently large $r$, the set $\{t\in [-\frac{r}{2},\frac{r}{2}] | f_t \in U  \}$ has measure at least $(1-\vn)r$.}

%\vspace{2mm}

 \vspace{2mm}
{\bf Theorem \ref{Main}}. {\em 
Let $U$ be any open neighborhood of $\g$ in the compact-open topology. Choose $t$ uniformly at random in $[-r,r]$.
The probability that $f_t$ is in $U$ tends to 1 as $r\ra \infty$.}
\vspace{2mm}

The motivation for this work comes from the study of the Kobayashi and Carath\'eodory metrics on Teichm\"uller space (see Section \ref{Teich}).
Let $\T$ denote the Teichm\"uller space of a finite-type orientable surface. 
A {\em Teichm\"uller disk} $\tau:\H \ra \T$ is a complex geodesic for the Kobayashi metric on $\T$.
It is an open problem to classify Teichm\"uller disks on which the
Kobayashi and Carath\'eodory metrics coincide.
To say that the metrics agree on $\tau(\H)$ means exactly that there is a {\em holomorphic retraction}
onto $\tau$, i.e. a holomorphic $\Psi:\T \ra \H$ so that $\Psi \circ \tau = \id_\H$. 

In recent work with Markovic \cite{GeMa}, we classify holomorphic retracts in the Teichm\"uller space of the five-times punctured sphere.
Key to our argument is the observation that certain Teichm\"uller disks $\tau$ factor as $$\H \stackrel{\Delta}{\ra} \H^n \stackrel{\mathcal{E}}{\ra} \T,$$ where $\Delta$
is the diagonal mapping and $\mathcal{E}$ is a particular naturally-defined holomorphic embedding.
If $\Psi:\T \ra \H$ is a holomorphic retraction onto $\tau(\H)$, then $f=\Psi \circ \mathcal{E}:\H^n \ra \H$ is a holomorphic retraction onto the diagonal,
i.e. $f$ is in $\mathcal{D}$. In \cite{GeMa}, we use the properties of $\mathcal{D}$ developed in this paper to glean information about holomorphic maps out of Teichm\"uller space.

The translation flow \eqref{tflow} should be viewed in the context of unipotent dynamics. 
The translation flow on $\mathcal{D}$ extends to an action of $\Aut(\H) \cong \PSL_2(\R)$ (See Section \ref{PAR}).
Equation \eqref{tflow} gives the action of the unipotent subgroup $$U=\lbr \lp \begin{array}{cc} 1& t\\ 0 &1 \end{array} \rp  \middle|  t\in \R \rbr.$$
Analogously, there is a natural $\PSL_2(\R)$ action on the unit cotangent bundle $T^*_1\T$ of Teichm\"uller space. The restriction of this action to $U$ is called the {\em horocycle flow}.
Our methods in \cite{GeMa} are summarized as follows: First, use results on horocycle flow in $T^*_1\T$ \cite{Sm} to reduce to an appropriate class of Teichm\"uller disks.
Next, use translation flow in $\mathcal{D}$ and the results of this paper to analyze retractions onto disks in that class.

Generalizing from the case of translations acting on holomorphic maps $\H^n \ra \H$, it is natural to ask the following equation:
Given two Hermitian symmetric spaces $X_1$ and $X_2$, what can one say about the dynamics of subgroups of $\Aut(X_1)\times \Aut(X_2)$ acting on subsets of the space of holomorphic functions $\mathcal{O}(X_1, X_2)$? To our knowledge, there is no previous work in the literature explicitly addressing this question. There has however been much interest in the dynamics of linear operators acting on holomorphic function spaces (see \cite{Ba}). In Section \ref{Ck} we use results \cite{Bo}\cite{Go} on linear dynamics to study the analogue of translation flow for maps $\C^k \ra \C$.
In this context, the flow is chaotic and behaves quite differently than the flow on maps $\H^k\ra \H$. 

The key tool in the proof of our main result is a multivariate version of the Schwarz lemma (see Section \ref{Schwarz}).
Our methods are inspired by Knese's work \cite{Kn} on extremal maps $\D^n \ra \D$.

\subsection{The Carath{\'e}odory and Kobayashi metrics on Teichm\"uller space}\label{Teich}

The Carath{\'e}odory pseudometric $d_C$ on a complex manifold $X$ assigns to two points $p,q\in X$ the distance
$$d_C(p,q) \equiv \sup_{f} d_{\H}(f(p), f(q)),$$ where the supremum is taken over all holomorphic maps
$f:X \ra \H$, and $d_\H$ is the Poincar\'e metric. In other words, $d_C$ is the smallest pseudometric on $X$ so that every holomorphic map from $X$ to $\H$
is length-decreasing.

The Kobayashi pseudometric $d_K$ on $X$ is defined in terms of maps $\H \ra X$. It is the largest pseudometric on $X$ so that every holomorphic map from $\H$ to $X$
is length-decreasing.

The Kobayashi and Carath\'{e}odory metrics on $\H^n$ are both given by
$$d_{\H^n}(z,w) = \max_{j} d_{\H}(z_j,w_j).$$  
In general, the Schwarz lemma implies $d_C\leq d_K$ for any complex manifold. However, it is usually difficult to determine if $d_C=d_K$ for a given complex manifold $X$.

In \cite{Ma}, Markovic proves that $d_C$ and $d_K$ do not agree on the Teichm\"{u}ller space of a closed orientable surface of genus $\geq 2$.
Let $\T$ be the Teichm\"uller space of a finite-type orientable surface.
Given a rational Strebel differential $\phi$ with characteristic annuli $\Pi_1,\ldots, \Pi_n$,
Markovic defines a holomorphic map $\mathcal{E}^{\phi}:\H^n \ra \T$.
The marked surface $\mathcal{E}^{\phi}(z_1,\ldots z_n)$ is constructed by applying the affine transformation $x+iy\mapsto x+z_j y$ to $\Pi_j$.
In particular, the restriction of $\mathcal{E}^\phi$ to the diagonal is the Teichm\"uller disk generated by $\phi$.  
Let $\alpha_j = \lp \int_{\Pi_j}\lpi \phi \rpi \rp / \|\phi\|_1 .$
Markovic proves the following:
\begin{Proposition}\label{Ma2}
 If the metrics $d_C$ and $d_K$ agree on the Teichm\"uller disk generated by $\phi$, then there is a holomorphic function $\Psi:\T\ra \H$ and a real constant
$T$ so that $f = \Psi\circ \mathcal{E}^{\phi}$ satisfies 
\begin{equation}\label{i}
f(\lambda,\ldots, \lambda) = \lambda \tag{A}
\end{equation} 
\begin{equation}\label{ii}
\pd{f}{z_j}(\lambda, \ldots, \lambda) = \alpha_j \tag{B}
\end{equation}
\begin{equation}\label{iii}
f(z_1+T, z_2+T,\ldots, z_n+T) = f(z_1, z_2,\ldots, z_n)+T \tag{C}
\end{equation}
for all $\lambda \in \H$, $(z_1,\ldots, z_n)\in \H^n$, $j=1,\ldots,n$.
\end{Proposition}

Markovic then proves
\begin{Proposition}\label{Ma}
For $n=2$, the only holomorphic $f:\H^2\ra \H$ satisfying conditions \eqref{i},\eqref{ii},\eqref{iii} is $f(z_1,z_2)=\alpha_1 z_1 + \alpha_2 z_2$.
\end{Proposition}
So if $\phi$ has exactly two characteristic annuli, there is a $\Psi: \T \ra \H$ such that $\Psi\circ \mathcal{E}^{\phi} = \alpha_1z_1+\alpha_2z_2$.
This criterion is then used to show that $d_C$ and $d_K$ do not agree on the Teichm\"uller disk generated by an $L$-shaped pillowcase
with rational edge lengths. 

As a corollary of our main result Theorem \ref{Main}, we obtain the generalization of Proposition \ref{Ma} to arbitrary $n$:

{\bf Corollary \ref{Cor}.}
{\em The only holomorphic $f:\H^n \ra \H$ satisfying \eqref{i},\eqref{ii},\eqref{iii} is} 
$f(z_1,\ldots, z_n) = \sum_{j=1}^n \alpha_j z_j$.

Taken together, Proposition \ref{Ma2} and Corollary \ref{Cor}
yield the following criterion for determining whether $d_C$ and $d_K$ agree on the Teichm\"uller disk generated by a rational Strebel differential. 
\begin{Proposition}\label{App}
Let $\phi$ be a rational Jenkins-Strebel differential, with characteristic annuli $\Pi_1,\ldots, \Pi_n$.
Suppose $d_C$ and $d_K$ agree on the Teichm\"uller disk generated by $\phi$.
Then there exists a holomorphic map $\Phi:\T \ra \H^n$
such that $$\Phi\circ \mathcal{E}^{\phi}(z_1,\ldots, z_n) = \alpha_1z_1+\cdots+\alpha_n z_n,$$
where $\alpha_j = \lp \int_{\Pi_j}\lpi \phi \rpi \rp /\|\phi\|_1 .$
\end{Proposition}

{\em Remark}:
Markovic showed that there are Teichm\"uller disks on which $d_C\neq d_K$. On the other hand, Kra \cite{Kr} proved that $d_C=d_K$ on every Teichm\"uller disk generated by
a holomorphic quadratic differential with no odd-order zeros.  
This raises a natural question: For which quadratic differentials do the Carath\'eodory and Kobayashi metrics on the corresponding disk agree?
A natural conjecture is that the converse of Kra's result holds: $d_C=d_K$ on a Teichm\"uller disk 
if and only if the generating differential has no odd-order zeros. In a recent paper \cite{GeMa} we prove this conjecture in the case of the five-times punctured sphere
and twice-punctured torus. Key to the proof is the fact that Proposition \ref{App} continues to hold without the rationality assumption.  
This fact in turn hinges on the main result Theorem \ref{Main} of this paper. (The weaker result Corollary \ref{Cor} is insufficient to deal with the irrational case.)   

\subsection{The Schwarz lemma and extremal maps}\label{Schwarz} 
Let $\D$ be the open unit disk in the complex plane.
The classical Schwarz lemma states that, if $f:\D \ra \D$ is holomorphic, then
\begin{equation}\label{CS}
 (1-\lpi z \rpi)^2 \lpi f'(z) \rpi \leq 1-\lpi f(z) \rpi^2,
\end{equation}
for all $z\in \D$. If equality holds in \eqref{CS} for some $z\in \D$, then it holds for all $z$. In this case, $f$ is a conformal automorphism of $\D$.

The Schwarz lemma has the following generalization for holomorphic maps $f$ from the polydisk
$\D^n = \D\times \cdots \times \D$ to $\D$ (see page 179 of \cite{Ru}):
\begin{equation}\label{S1}
\sum_{j=1}^n (1-\lpi z_j \rpi ^2)\lpi  \pd{f}{z_j}(z) \rpi \leq 1-\lpi f(z) \rpi^2,
\end{equation}
for every $z=(z_1,\ldots, z_n)\in \D^n$.
To understand $\eqref{S1}$, we recall the following definitions: A {\em balanced disk} in $\D^n$ is a copy of
$\D$ embedded in $\D^n$ by a map of the form
$$\Phi: z\mapsto \lp \phi_1(z),\ldots, \phi_n(z)\rp,$$
where $\phi_i \in \Aut(\D)$. A balanced disk $\Phi$ is called {\em extreme} for $f$ if the restriction $f\circ \Phi$ is in $\Aut(\D)$.   
The content of $\eqref{S1}$ is that the restriction of $f$ to every balanced disk satisfies the classical Schwarz lemma.
Equality in $\eqref{S1}$ means that $z$ is contained in some extreme disk for $f$.  

The {\em extreme set} $X(f)$ is the union of the extreme disks of $f$. In other words, $X(f)$ is the set of points $z\in \D^n$ for which equality holds
in \eqref{S1}. In \cite{Kn}, Knese classifies maps $f:\D^n \ra \D$ for which $X(f)=\D^n$. Such maps are called {\em everywhere extremal}, or simply {\em extremal}.
Knese shows that extremal maps $\D^n\ra \D$ form a special class of rational functions parameterized by $(n+1)\times (n+1)$ symmetric unitary matrices.

The upper half-plane $\H$ is conformally equivalent to $\D$ via the Cayley transform $z\mapsto \frac{i-z}{i+z}$.  
For holomorphic maps $f:\H^n \ra \H$, the generalized Schwarz lemma becomes
\begin{equation}\label{Schwarz}
\sum_{j=1}^n \im(z_j)\lpi  \pd{f}{z_j}(z) \rpi \leq \im f(z).
\end{equation}

\subsection{The families $\mathcal{D},\mathcal{C}$}\label{D}
Consider the family $\mathcal{D}$ of holomorphic maps $f:\H^n \ra \H$ which restrict to the identity on the diagonal: 
\begin{equation} \label{diag}
f(\lambda ,\ldots, \lambda)=\lambda
\end{equation}
for all $\lambda\in \H$. $\mathcal{D}$ is a natural class to consider; it is the collection of maps $\H^n\ra \H$ with a distinguished extreme disk.
After pre- and post-composing by biholomorphisms, any holomorphic map $\H^n\ra \H$ with an extreme disk becomes an element of $\mathcal{D}$. 

Differentiating both sides of \eqref{diag} with respect to $\lambda$ yields $$\sum_{j=1}^n\pd{f}{z_j}(\lambda,\ldots, \lambda) = 1.$$
But by the generalized Schwarz lemma \eqref{Schwarz},
$$\sum_{j=1}^n\lpi \pd{f}{z_j}(\lambda,\ldots, \lambda)\rpi \leq 1.$$
So $\pd{f}{z_j}(\lambda,\ldots, \lambda)\geq 0$
for all $\lambda\in \H$ and $j=1,\ldots, n$. By the open mapping theorem, $\lambda \mapsto \pd{f}{z_j}(\lambda,\ldots, \lambda)$ is constant.
So $f$ satisfies $$\pd{f}{z_j}(\lambda,\ldots, \lambda) = \alpha_j$$ for all $\lambda\in \H$, for some collection of nonnegative constants $\alpha_j$ summing to
1. 

In the rest of the paper, we assume without loss of generality that $\alpha_j = \frac{1}{n}$. To reduce the general case to this one, suppose $f\in \mathcal{D}$
and $\pd{f}{z_j}(i,\ldots, i) = \alpha_j$. Define $g\in\mathcal{D}$ by $$g(z) = \sum_{j=1}^{n} \lp \frac{1-\alpha_j}{n-1}\rp z_j.$$
Then $$\tilde{f} = \frac{1}{n}f + \frac{n-1}{n}g$$ is in $\mathcal{D}$ and satisfies $\pd{\tilde{f}}{z_j}(i,\ldots, i) = \frac{1}{n}$.
Since $g$ is invariant under the translation flow, it suffices to consider the translation orbit of $\tilde{f}$.

With these considerations in mind, we define $\mathcal{C}$ to be the family of holomorphic maps $\H^n \ra \H$ satisfying
\begin{equation}\label{A}
f(\lambda,\ldots, \lambda) = \lambda, \tag{A}
\end{equation}
\begin{equation}\label{B}
\pd{f}{z_j}(\lambda,\ldots, \lambda) = \frac{1}{n}, \tag{B}
\end{equation}
for all $\lambda \in \H$ and $j=1,\ldots, n.$

When convenient, we view $\mathcal{C}$ as the family of maps $\D^n \ra \D$ satisfying the same conditions.
(Conjugation by the Cayley transform $\H\ra \D$ preserves $\eqref{A},\eqref{B}$.)

{\em Remark}: Conditions $\eqref{A}$ and $\eqref{B}$ hold for all $\lambda\in \H$ iff they both hold for some $\lambda \in \H$.

\subsection{Extremal maps in dimension two}\label{Exa}
In \cite{Kn}, Knese showed that extremal maps $g:\D^2 \ra \D$ satisfying $g(0,0)=0$
are all of form 
$$
g(z,w)= \mu\frac{az+bw-zw}{1-\ol{b}z-\ol{a}w},
$$
where $|\mu| = |a|+|b|=1$.
Imposing $f(\lambda,\lambda)=\lambda$ and $\pd{f}{z}(\lambda,\lambda)=\pd{f}{w}(\lambda,\lambda) = \frac{1}{2}$, we find that the extremal elements of $\mathcal{C}$
are the functions of form
$$g_\nu(z,w) = \frac{\nu(\frac{z}{2}+\frac{w}{2})-zw}{\nu-(\frac{z}{2}+\frac{w}{2})} $$
with $\nu \in \partial \D$.

A direct computation shows that, for any $\gamma\in \Aut(\D)$, $$\gamma\cdot g_\nu = g_{\gamma(\nu)},$$
where $(\gamma\cdot g_\nu)(z_1,z_2) = \gamma g_\nu (\gamma^{-1}z_1, \gamma^{-1}z_2).$
Thus, the set of extremals in $\mathcal{C}$ is in $\Aut(\D)-$equivariant bijection with $\partial \D$.
%In particular, $g_\nu$ is invariant under the conjugation action of $$\stab_\nu =\{\gamma\in \Aut(\D) | \gamma(\nu) = \nu \}.$$

{\em Remark}: The situation for $n>2$ is more complicated; one can show using Knese's classification of extremals that the extremals in $\mathcal{C}$ constitute a manifold of dimension
$\frac{n(n-1)}{2}$.

Conjugating by the Cayley transform, we get a description of the extremal maps $\H^2\ra \H$ in $\mathcal{C}$.
They are the functions of form
$$h_r(z,w) = \frac{r(\frac{z}{2}+\frac{w}{2})-zw}{r-(\frac{z}{2}+\frac{w}{2})},$$
with $r\in \partial \H=\R\cup \{\infty\}$. In particular,
$$h_\infty(z,w) = \frac{z}{2} + \frac{w}{2}.$$

One can check that the extreme disks for $h_\infty$ are precisely those of form $\{\lp z,az+b\rp |z\in\H\}$,
where $a>0$ and $b\in\R$.
It follows, more generally, that the extreme disks for $h_r$ are those of form $\{\lp z,\phi(z)\rp| z\in \H\}$,
with $\phi\in \stab(r)$.

{\bf Example 1.3}
In \cite{Kn}, Knese constructed a holomorphic map $\D^2 \ra \D$ which has two extreme disks, yet is not everywhere extremal.
Below, we give an example of a map $\H^2 \ra \H$ which is extremal on every disk of the form $\{(z,a z)| z\in \H\}$ with $a>0$,
yet is not everywhere extremal.

Given $r,s\in \partial \H$, $\stab(r)\cap \stab(s)$ is the set of isometries preserving the hyperbolic geodesic with endpoints $r,s$.
For example, $\stab(0)\cap \stab(\infty)$ consists of isometries preserving the positive imaginary axis; these are of form $z\mapsto az$ with $a>0$.
So the disks $D_a = \{(z, az)| z\in \H\}$ are extreme for both $h_\infty(z,w)=\frac{z+w}{2}$ and $h_0(z,w)=\frac{2zw}{z+w}$. 
In fact, the $D_a$ are extreme for any convex combination
$$f^t=t h_\infty + (1-t)h_0,$$
with $t\in(0,1)$. 
Indeed,
$$f^t(z,az)=\lp t\frac{1+a}{2}+ (1-t)\frac{2a^2}{1+a}\rp z.$$
So the extreme set $X(f^t)$ contains a set of real dimension 3. Yet $f^t$ is not everywhere extremal, as $f^t\neq h_r$ for any $r\in \partial \H$.

\subsection{Translation flow in dimension 2}
In dimension 2, $\mathcal{C}$ can be parameterized explicitly using Nevanlinna-Pick interpolation on the bidisk.
The maps $\D^2 \ra \D$ belonging to $\mathcal{C}$ are precisely those of form
\begin{equation}\label{AM}
f(z,w) = \frac{1}{2}(z+w)+\frac{1}{4}(z-w)^2 \frac{\Theta(z,w)}{1-\frac{1}{2}(z+w)\Theta(z,w)},
\end{equation}
where $\Theta$ is any holomorphic map from $\D^2$ to the closed disk $\ol{\D}$. 
(See page 189 of \cite{Ag}.)

To parameterize maps $\H^2\ra \H$ in $\mathcal{C}$, we conjugate \eqref{AM} by the Cayley transform.
We get the same general form, with $\Theta$ any holomorphic map from $\H^2$ to the closure $\ol{\H}$ of $\H$ in the Riemann sphere.
Substituting $\Theta = -\frac{1}{\Phi}$, \eqref{AM} becomes
\begin{equation}\label{Be}
f(z,w) = \frac{\frac{z+w}{2}\cdot \Phi(z,w)+zw}{\Phi(z,w) + \frac{z+w}{2}},
\end{equation}
The extremal map $h_r$ corresponds to $\Phi \equiv -r$. In particular, $h_\infty(z,w)=\frac{z+w}{2}$ correponds to $\Phi \equiv \infty$.

Applying translation flow to \eqref{Be} yields
\begin{equation}\label{Bet}
f(z-t,w-t) + t = \frac{\frac{z+w}{2}\cdot \lbrac \Phi(z-t,w-t)-t\rbrac +zw}{\lbrac \Phi(z-t,w-t)-t\rbrac + \frac{z+w}{2}}.
\end{equation}

One can show that for randomly chosen real $t$, $\lpi \Phi(z-t, w-t)-t \rpi$ is very large, so that \eqref{Bet} is very close to $\frac{z+w}{2}$. 
This yields a proof of Theorem \ref{Main}
in dimension 2. 

\subsection{Translation flow for maps $\C^k \ra \C$}\label{Ck}
Let $\mathcal{D}'$ denote the space of holomorphic maps $f:\C^k \ra \C$ satisfying $f(z,\ldots, z)=z$ for all $z\in \C$.
Define translation flow on $\mathcal{D}'$ by the same formula
$$
f_t(z_1,\ldots, z_k) = f(z_1-t, \ldots, z_k-t)+t
$$
as the flow on $\mathcal{D}$. 

The main results of this paper state that translation flow on $\mathcal{D}$ is ``unchaotic."
Theorem \ref{Main} asserts that the orbit any $f\in \mathcal{D}$ is concentrated at a single point,
while \ref{Cor} states that the periodic points lie in a finite-dimensional subspace of $\mathcal{D}$.
In stark contrast, the flow on $\mathcal{D}'$ has orbits which equidistribute; moreover,
the set of periodic points is dense. This contrast should be viewed
in light of the fact that, unlike $\mathcal{D}$, the space $\mathcal{D}'$ is not compact.

\begin{Proposition}\label{Dense1}
There is a probability measure $\mu$ on $\mathcal{D}'$ which is ergodic with respect to translation flow
and whose support is the entire space $\mathcal{D}'$.
In particular, a dense set of $f\in \mathcal{D}'$ have $\mu$-equidistributed orbits under translation flow.
\end{Proposition}

{\em Remark:}
Another way of stating the main result Theorem \ref{Main} is that any ergodic probability measure for translation flow on $\mathcal{D}$ is a
delta measure supported at a point $g\in \mathcal{D}$ of form
$g(\mathbf{z}) = \sum_j \alpha_j z_j$ (see Proposition \ref{ERG}).

\begin{Proposition}\label{Dense2}
The set of periodic points for the translation flow on $\mathcal{D}'$ is dense.
\end{Proposition}

Propositions \ref{Dense1} and \ref{Dense2} follow easily from the following results on linear dynamics:

\begin{Proposition}[{\bf Bonilla, Grosse-Erdmann \cite{Bo}}]\label{P1}
Let $L$ be any continuous linear operator on $\mathcal{O}(\C^n)$ which commutes
with the differential operators $\pd{}{z_1},\ldots, \pd{}{z_n}$. 
Then $L$ is ergodic with respect to a full-support probability measure. 
\end{Proposition}  

\begin{Proposition}[{\bf Godefrey, Shapiro \cite{Go}}]  \label{P2}
Under the hypotheses of Proposition \ref{P1}, $L$ has a dense set of periodic points. 
\end{Proposition}

{\em Proof of Propositions \ref{Dense1} and \ref{Dense2}:}
Let $S:\mathcal{D}' \ra \mathcal{D}'$ be the time-one translation $f\mapsto f_1$.
It suffices to show that $S$ has a dense set of periodic points and an ergodic probability measure
$\mu$ with full support. (To obtain the desired flow-invariant measure, average $\mu$ over the flow from time 0 to time 1.)  

Propositions \ref{P1},\ref{P2} apply to the operator $L$ on $\mathcal{O}(\C^n)$ defined by
$$
L\phi(z_1,\ldots, z_n) = \phi(z_1-1,\ldots, z_n-1).
$$

It thus suffices to exhibit a continuous surjection
$$
\mathcal{O}(\C^n) \ra \mathcal{D}'.
$$
intertwining the actions of $L$ and $S$.

To this end, define $g\in \mathcal{O}(\C^k)$ by $g(z_1,\ldots, z_n) = \frac{1}{n}\sum_j z_j$,
and let $\mathbf{1}\in \C^n$ denote the vector with all entries equal to 1.
The map $F:\mathcal{O}(\C^n) \ra \mathcal{D}'$ associating to each $\phi \in \mathcal{O}(\C^n)$ the function 
$$f(\z) = \phi(\z) - \phi(g(\z){\mathbf 1}) + g(\z)$$
is the desired surjection. It is easy to check $F$ intertwines the actions of $L$ and $S$.
Moreover, the map
$\mathcal{D}'\ra \mathcal{O}(\C^n)$ sending $f\in \mathcal{D}'$ to
$\phi(\z) = f(\z)-g(\z)$ is a right inverse for $F$.
\qed 
   
\subsection{Outline}
The rest of the paper will focus on the proof of our main result, Theorem \ref{Main}.
The key observation is that $\g(z) = \frac{1}{n}\sum_{j=1}^n z_j$ is an everywhere extremal map from $\H^n$ to $\H$.

In Section \ref{Ext}, we show that extremals in $\mathcal{C}$
are extreme points of $\mathcal{C}$, in the sense of convex analysis. More precisely, we prove 

{\bf Proposition \ref{Ex2}.}
{\em If $g\in \mathcal{C}$ is extremal and $\mu$
is a Borel probability measure on $\mathcal{C}$ such that $$\int_{\mathcal{C}} f(z) d\mu(f) = g(z)~ \forall z \in \H^n,$$ then $\mu$ is the Dirac measure $\delta_g$ concentrated at the point $g\in \mathcal{C}$.}

Then, in Section \ref{Trans} we show that the average of any $f\in \mathcal{C}$ over the translation flow is $\g(z)=\frac{1}{n}\sum_{j=1}^n z_j$. That is, we prove

{\bf Proposition \ref{AV}.}
{\em Let $f\in \mathcal{C}$. For each $t\in \R$, define $f_t(z_1,\ldots, z_n) = f(z_1-t,\ldots, z_n-t)+t$. 
Then $\frac{1}{2r}\int_{-r}^{r} f_t(z) dt$ converges locally uniformly to $\g(z)$ as $r\ra\infty$.
}

In Section \ref{MR}, we prove the main result.
To apply Proposition \ref{Ex2}, we consider the measure $\mu_r$ on $\mathcal{C}$ obtained by pushing forward the uniform probability measure on $[-r,r]$ via the map $t\mapsto f_t$.
The desired result is that $\mu_r\ra \delta_\g$ as $r\ra\infty$. Propositions \ref{Ex2}, \ref{AV} imply that $\delta_\g$
is the only accumulation point of $\{ \mu_r\}_{r>0}$. The main result then follows by the Banach-Alaoglu theorem.

In Section \ref{PAR}, we rephrase our results in a more invariant form, in terms of the conjugation action of 
$\PSL_2(\R)$ on $\mathcal{D}$. In Section \ref{Six}, we establish a rigidity result used in the proof of Proposition \ref{AV}, and in the Appendix,
we discuss generalizations of the classical polarization principle. 
 
\section{Convexity and extreme points}\label{Ext}

Let $\mathcal{C}$ be the family of holomorphic maps $\H^n \ra \H$ satisfying
\begin{equation}\label{A}
f(\lambda,\ldots, \lambda) = \lambda, \tag{A}
\end{equation}
\begin{equation}\label{B}
\pd{f}{z_j}(\lambda,\ldots, \lambda) = \frac{1}{n}, \tag{B}
\end{equation}
for all $\lambda \in \H$ and $j=1,\ldots, n.$

Recall that an extremal map $g:\H^n \ra \H$ is a holomorphic function satisfying
$$
\sum_{j=1}^n \im(z_j)\lpi  \pd{g}{z_j}(z) \rpi = \im g(z)
$$
for all $z=(z_1,\ldots, z_n)\in \H^n$.

Observe that $\mathcal{C}$ is a convex subset of the holomorphic functions on $\H^n$.
Our next result is that every extremal in $\mathcal{C}$ is an extreme point in the sense of convex analysis. 

\begin{Proposition}\label{Ex}
Suppose $g\in \mathcal{C}$ is extremal. If $g = tf_1 + (1-t)f_2$, with $f_i \in \mathcal{C}$ and $t\in (0,1)$,
then $f_1=f_2=g$.
\end{Proposition} 

\noindent {\em Proof}: 
We have
\begin{align}
t\> \im (f_1) + (1-t)\im (f_2) &= \sum_{j=1}^n \im (z_j)   \lpi t\pd{f_1}{z_j}+(1-t)\pd{f_2}{z_j} \rpi \label{Eq}\\
                                         &\leq \sum_{j=1}^n \im (z_j) \lbrac   t \lpi \pd{f_1}{z_j} \rpi + (1-t)\lpi \pd{f_2}{z_j} \rpi \rbrac \nonumber \\
                                         &\leq t\> \im (f_1) + (1-t)\im (f_2),\nonumber
\end{align}
where in the first line we've used that $g$ is extremal, 
and in the third we've applied \eqref{Schwarz} to $f_1,f_2$.
Thus,
$$\lpi t\pd{f_1}{z_j}(z)+(1-t)\pd{f_2}{z_j}(z) \rpi =    t \lpi \pd{f_1}{z_j}(z) \rpi + (1-t)\lpi \pd{f_2}{z_j}(z) \rpi $$
for $j=1,\ldots, n$ and all $z\in \H^n$. 

So $$\lp \pd{f_1}{z_j}\rp \lp \pd{g}{z_j}\rp^{-1}\geq 0,$$
whenever $\pd{g}{z_j}\neq 0$, and similarly for $f_2$. 
Let $U\subset \H^n$ be the complement of the zero set of $\pd{g}{z_j}$.
By $\eqref{B}$, $\pd{g}{z_j}$
is not identically zero, so $U$ is a dense connected subset of $\H^n$.
The open mapping theorem now implies that $\lp \pd{f_1}{z_j}\rp \lp \pd{g}{z_j}\rp^{-1}$ is a nonnegative constant on $U$.
Again by $\eqref{B}$,
$$ \pd{f_1}{z_j}= \pd{g}{z_j}$$
on $U$ and, thus, on all of $\H^n$.
Since the first derivatives of $f_1$ and $g$ are the same, $f_1$ and $g$ differ by a constant.
By $\eqref{A}$, $f_1=g$. Similarly, $f_2 = g$. \qed

The last result implies that if a finite convex combination $$g=\sum_k t_k f_k$$ of elements of $\mathcal{C}$
is extremal, then the $f_k$ are all equal to $g$.
We will show, more generally, that if $\mu$ is a Borel probability measure on $\mathcal{C}$
such that $$g=\int_\mathcal{C} f d\mu(f)$$ is extremal, then $\mu = \delta_g$.
Before we consider Borel measures on the space $\mathcal{C}$, we need to understand the space's basic topological properties.
\begin{Proposition}\label{Compact}
The family $\mathcal{C}$ is compact and metrizable in the compact-open topology.
\end{Proposition}

\noindent{\em Proof}:
Metrizability is standard: Choose a compact exhaustion $K_1,K_2,\ldots$ of $\H^n$, and set $d_j(f,g) = \sup_{z\in K_j} |f(z)-g(z)|$.
Then the metric
$$d(f,g) = \sum_{j=1}^\infty 2^{-j} \frac{d_j(f,g)}{1+d_j(f,g)}$$
induces the compact-open topology.

To prove compactness, we need to show that $\mathcal{C}$ is precompact and closed in $\mathcal{O}(\H^n)$.
By the definition of the Carath{\'e}odory metric, any holomorphic map $\H^n\ra \H$ decreases Carath{\'e}odory distance.
Thus, every $f \in \mathcal{C}$ satisfies
$$d_{\H}(f(z_1,\ldots, z_n), i) \leq d_{\H^n} \lp  (z_1,\ldots, z_n), (i,\ldots, i)  \rp.$$
The right side of the inequality is continuous in the $z_j$. So $\mathcal{C}$ is locally uniformly bounded and thus precompact.
The inequality also implies that any accumulation point of $\mathcal{C}$ has image contained in $\H$. Furthermore, $\eqref{A}$ and $\eqref{B}$ are closed conditions.
Thus, $\mathcal{C}$ is closed in $\mathcal{O}(\H^n)$. 
 \qed

\vspace{2mm}

Let $\mu$ be a Borel probability measure on $\mathcal{C}$. For each $z\in \H^n$,
the evaluation map $f\mapsto f(z)$ is a continuous function on the compact space $\mathcal{C}$.
So the evaluation map is $\mu$-integrable. We denote its integral by $\int_{\mathcal{C}}f(z)d\mu(f)$.

\begin{Proposition}\label{Ex2}
Suppose $g\in \mathcal{C}$ is extremal. Let $\mu$ be a Borel probability measure on $\mathcal{C}$. Suppose
$\int_{\mathcal{C}} f(z) d\mu(f)=g(z)$ for all $z\in \mathcal{\H}^n$. Then $\mu$ is $\delta_g$, the Dirac measure
concentrated at $g$.
\end{Proposition}

\noindent{\em Proof}:
Though this result can be derived as a formal consequence of Proposition \ref{Ex}, we prefer to give a direct proof.

The proof is similar to that of Proposition \ref{Ex}. To establish the analog of equality \eqref{Eq}, we need to differentiate
$\int_{\mathcal{C}} f(z) d\mu(f)$ under the integral sign; Proposition \ref{Compact} implies
that the family $\{\pd{f}{z_j} | f\in \mathcal{C} \}$ is locally uniformly bounded, which justifies switching $\int$ and $\pd{}{z_j}$.

Let $U$ be the complement of the zero set of $\pd{g}{z_j}$. Fix $z\in U$. Arguing as before, we get
\begin{equation}\label{RAT}
\lp \pd{f}{z_j}(z)\rp \lp \pd{g}{z_j}(z)\rp^{-1}\geq 0,
\end{equation}
for $\mu$-almost-every $f$. A countable intersection of full-measure subsets of $\mathcal{C}$ has full measure.
Thus, for $\mu$-a.e. $f$, \eqref{RAT} holds at all $z\in U$ with rational coordinates.
By continuity, $\mu$-a.e. $f$ satisfies \eqref{RAT} on $U$. We conclude that $\mu$-a.e. $f$ is equal to $g$.
This means that $\mu=\delta_g$. \qed

\section{Averaging over translations}\label{Trans}

Let $f:\H^n \ra \H$ be a holomorphic map.
For each $t \in \R$, we define
$$f_t(z_1,\ldots, z_n) = f(z_1-t,\ldots, z_n-t) + t.$$
The action $(f,t)\mapsto f_t$ is the translation flow on $\mathcal{O}(\H^n)$.
The family $\mathcal{C}$ is invariant under the translation flow.

For each $f\in\mathcal{C}$ and $r>0$, we define the average
$\mathcal{A}_r[f]\in \mathcal{C}$ by
$$\mathcal{A}_r[f](z)  = \frac{1}{2r}\int_{-r}^{r} f_t(z)dt.$$
One might expect that averaging $f\in \mathcal{C}$ over the entire flow yields an invariant element.
This is indeed the case:   
\begin{Proposition}\label{AV}
For each $f\in \mathcal{C}$, $\mathcal{A}_{r}[f]$ converges locally uniformly to
$\mathbf{g}(z)= \frac{1}{n}\sum_{j=1}^n z_j$ as $r\ra\infty$.
\end{Proposition}

\noindent{\em Proof:}
Fix $z\in \H^n$.
By Proposition \ref{Compact}, there is a $C(z)>0$ so that 
\begin{equation}\label{Bound}
\lpi f(z) \rpi < C(z)
\end{equation}
for all $f\in \mathcal{C}$.

Fix $s\in \R$. 
We use \eqref{Bound} to compare $\mathcal{A}_r[f]$ and the translate $\lp \mathcal{A}_r[f] \rp_s$:
\begin{align*}
\lpi \mathcal{A}_r[f](z)- \lp \mathcal{A}_r[f]\rp_s(z) \rpi &= \frac{1}{2r} \lpi \int_{-r}^{-r+s} f_t(z)dt - \int_{r}^{r+s} f_t(z)dt \rpi \\
												 &\leq \frac{s}{r}C(z).
\end{align*} 
Thus, any limit point of the family $\{\mathcal{A}_r[f]\}_{r>0}$ along a sequence with $r \ra \infty$ is invariant under all translations.
But, as we will show in Proposition \ref{AF}, the only translation-invariant element of $\mathcal{C}$ is $\mathbf{g}$.
Since $\mathcal{C}$ is sequentially compact, we get the desired result.\qed

\section{The main result}\label{MR}

We now use Propositions \ref{Ex2}, \ref{AV} and the Banach-Alaoglu theorem to prove the main result.

\begin{Theorem}\label{Main}
Suppose $f:\H^n \ra \H$ is holomorphic and satisfies $f(\lambda,\ldots, \lambda)=\lambda$ for all $\lambda\in \H$. Let $\alpha_j = \pd{f}{z_j}(i,\ldots, i)$, and
define $\g(z) = \sum_{j=1}^n \alpha_j z_j$.
Fix $\vn>0$, and let $U$ be any open neighborhood of $\g$ in the compact-open topology. Then for sufficiently large $r$, the set $\{t\in [-\frac{r}{2},\frac{r}{2}] | f_t \in U  \}$ has measure at least $(1-\vn)r$.
\end{Theorem}

\noindent{\em Proof:}
We may assume without loss of generality that $\alpha_j = \frac{1}{n}$ for $j=1,\ldots, n$. So $f\in \mathcal{C}$.
Let $\mu_r$ be the pushforward to $\mathcal{C}$ of the uniform probability measure on $[-r,r]$, via the continuous map $t\mapsto f_t$.
Then the desired result is equivalent to the assertion that  $\mu_r \ra \delta_{\mathbf{g}}$ weakly as $r\ra \infty$.

By the Banach-Alaoglu theorem, the space of Borel probability measures on the compact metric space $\mathcal{C}$ is sequentially compact.
It thus suffices to show that any limit point $\mu$ of $\{\mu_r\}_{r>0}$ along a sequence with $r\ra \infty$ is $\delta_\g$.
Proposition \ref{AV} says that $\int_\mathcal{C} h(z) d\mu_r(h) \ra \g(z)$, as $r\ra\infty$. 
So $\mu$ satisfies $$\int_\mathcal{C} h(z) d\mu(h) = \g(z)$$
for all $z\in \H^n$.
By Proposition \ref{Ex2}, $\mu = \delta_{\g}$. This completes the proof.\qed

\vspace{2mm}

The Birkhoff ergodic theorem yields following restatement of the main result.
\begin{Proposition}\label{ERG}
The only invariant measure for translation flow on $\mathcal{C}$ is the delta measure $\delta_\g$.
\end{Proposition}
{\em Remark}:
We do not know if $\lim_{t\ra \infty} f_t = \g$ for all $f\in \mathcal{C}$. 
  
As a corollary to the main result, we obtain the generalization of Proposition \ref{Ma} to maps $\H^n\ra \H$.
\begin{Corollary}\label{Cor}
Suppose $f:\H^n\ra \H$ is holomorphic and satisfies $f(\lambda,\ldots, \lambda)=\lambda$ for all $\lambda\in \H$. 
Suppose in addition that $f(z_1+T,\ldots, z_n+T) = f(z_1,\ldots, z_n)+T$ for some $T>0$ and all $(z_1,\ldots,z_n)\in \H^n$.
Then $f$ is equal to the function $\g(z) = \sum_{j=1}^n\alpha_j z_j $, where $\alpha_j = \pd{f}{z_j}(i,\ldots,i)$.
\end{Corollary}
\noindent{\em Proof:} Assume WLOG $\alpha_j=\frac{1}{n}$.
The hypothesis on $f$ means that it is a periodic point of the translation flow, with period $T$.
Thus, $\mu_{T} = \lim_{r\ra \infty}\mu_{r} = \delta_{\g}$.
Since $t\mapsto f_t$ is continuous, it follows that $f_t = \mathbf{g}$ for all $t\in [-T, T]$.
In particular, $f = \mathbf{g}$, as claimed.

\section{Unipotent subgroups acting on $\mathcal{D}$}\label{PAR}
In this section, we restate our results in terms of the action of $\Aut(\H)$ on $\mathcal{D}$.

The group $\Aut(\H)\cong \PSL_2(\R)$ acts on $\mathcal{D}$ by conjugation:
An element $\gamma \in \PSL_2(\R)$ sends $f\in\mathcal{D}$ to the function $\gamma\cdot f$ given by
$$
(\gamma \cdot f)(z_1,\ldots, z_n) = \gamma f(\gamma^{-1}z_1,\ldots, \gamma^{-1}z_n).
$$
By the chain rule,
$\gamma\cdot f$ has the same first partials at $(i,\ldots, i)$ as $f$.
So $\mathcal{C}$ is invariant under the action.

An element of $\PSL_2(\R)$ is called {\em unipotent} (or {\em parabolic}) if it fixes exactly one point in $\partial \H$.
A {\em unipotent subgroup} of $\PSL_2(\R)$ is a nontrivial one-parameter subgroup whose non-identity elements are unipotent.
Every unipotent subgroup is conjugate to the group of translations $z\mapsto z+t$. 

\vspace{2mm}

The following generalization of our results is immediate:
\begin{Theorem}\label{Par}
Let $\mathcal{D}$ be the family of holomorphic maps $\H^n\ra \H$ which restrict to the identity on the diagonal. 
Let $f\in \mathcal{D}$. For each $j$, $\lambda\mapsto \pd{f}{z_j}(\lambda,\ldots, \lambda)$ is identically equal to some nonnegative constant $\alpha_j$. \\
Let $\{\gamma_t\}\subset \PSL_2(\R)$ be a unipotent subgroup.
There is a unique $\gamma_1$-invariant holomorphic $\mathbf{g}\in \mathcal{D}$ satisfying 
$\pd{\mathbf{g}}{z_j}(\lambda,\ldots, \lambda) = \alpha_j$ for all $\lambda\in \H$ and $j=1,\ldots, n$.
\\Let $\mu_r$ be the pushforward to $\mathcal{D}$ of the uniform measure on $[-r,r]$, by the map $t\mapsto \gamma_t\cdot f$.
Then $\mu_r \ra \delta_{\mathbf{g}}$ weakly as $r\ra\infty$. 
\end{Theorem}

{\em Remark:}
Theorem \ref{Par} holds exactly as stated with $\H$ replaced by $\D$.

\section{A rigidity result}\label{Six}
Below, we establish the rigidity result we used in the proof of Proposition \ref{AV}, namely that any $f\in \mathcal{D}$ which is invariant under all translations is a convex combination of the coordinate functions.

\vspace{2mm}

First, we need a lemma.
\begin{Lemma}\label{Harmonic}
Let $\phi:\C \ra \R$ be a harmonic function with $\phi(0) = \pd{\phi}{x}(0) = \pd{\phi}{y}(0) = 0$. Suppose there is a $C>0$ so that 
 $\phi(z) \geq -C |z|$ for all $z\in \C$. Then $\phi$ is identically zero.
\end{Lemma}
\noindent {\em Proof:}
The idea is to use the Poisson integral formula to show that $\phi$ has sublinear growth.

Write $\phi = \phi_+ - \phi_{-},$ where $\phi_+(z) = \max\{0, \phi(z)\}$, and $\phi_-(z) = \max\{0, -\phi(z)\}$.
Fix $r>0$, and
set $$A = \int_{0}^{1} \phi_+(re^{2\pi i\theta})d\theta,~B = \int_{0}^{1} \phi_-(re^{2\pi i\theta})d\theta.$$
By the mean value property, $A-B = \phi(0)=0$. 
We compute
\begin{align*}
\int_0^1 |\phi(re^{2\pi i\theta})|d\theta &= A+B\\ &= 2B\\ &= 2\int_0^1\phi_-(re^{2\pi i\theta})d\theta\\ &\leq 2C r,
\end{align*}
where in the last inequality, we've used $\phi(z) \geq -C |z|$.
Now, for any $z$ with $|z| = \frac{r}{2},$ the Poisson integral formula for the ball $B_r(0)$ yields 
$$
\lpi \phi(z) \rpi = \lpi \int_{0}^1 \frac{r^2-\lp\frac{r}{2}\rp^2}{r\lpi z-re^{2\pi i \theta} \rpi} \phi(re^{2\pi i \theta}) d\theta \rpi 
\leq  \sup_{\theta\in [0,2\pi]} \lp \frac{3r}{4\lpi z-re^{2\pi i \theta} \rpi}  \rp \cdot \int_0^1 |\phi(re^{2\pi i\theta})|d\theta \leq 3C r.
$$ 

Since $r$ was arbitrary, we have $|\phi(z)| \leq 6C |z|$ for all $z$. Since $\phi$ is harmonic and has sublinear growth, $\phi$ is affine, that is, $\phi(x+iy) = ax+by+c$ for some $a,b,c\in \C$.
(Indeed, the higher derivatives of $\phi$ at 0 vanish, as we can see  
by differentiating Poisson's formula on $B_r(0)$ under the integral and letting $r$ tend to infinity.) 
By assumption, $\phi$ and its first derivatives vanish at the origin, so $\phi$ is identically 0.
   
\qed

We now prove the main result of this section.
\begin{Proposition}\label{AF}
Fix positive constants $\alpha_j$ with $\sum_{j=1}^n \alpha_j=1$.
Let $f:\mathcal{\H}^n \ra \H$ be a holomorphic function satisfying
\begin{equation}\label{a}
f(\lambda,\ldots, \lambda) = \lambda, \tag{A}
\end{equation}
for all $\lambda \in \H.$
\begin{equation}\label{b}
\pd{f}{z_j}(\lambda,\ldots, \lambda) = \alpha_j, \tag{B}
\end{equation}
for all $\lambda\in \H$ and $j=1,\ldots, n.$
\begin{equation}\label{c}
f(z_1+t,\ldots,z_n+t) = f(z_1,\ldots, z_n)+t, \tag{C}
\end{equation}
for all $(z_1,\ldots, z_n) \in \H^n$ and all $t\in \R$.
Then $f$ is the function $f(z_1,\ldots, z_n) = \sum_{j=1}^n \alpha_jz_j$.
\end{Proposition}

\noindent {\em Proof:}
As usual, we assume $\alpha_j=\frac{1}{n}$.
The idea is to first show that $f$ is of form $$\frac{1}{n}\sum_{j=1}^n z_j + H(z_2-z_1, z_3-z_2,\ldots, z_{n}-z_{n-1}),$$
for some holomorphic $H:\C^{n-1}\ra \C$. Then we use Lemma \ref{Harmonic} to show that $H\equiv 0$.

Let $$g(z_1,\ldots, z_n) = f(z_1,\ldots, z_n)-  \frac{1}{n}\sum_{j=1}^n z_j.$$ In terms of $g$, conditions \eqref{a}, \eqref{b}, \eqref{c} become 
\begin{equation}\label{A'}
g(\lambda,\ldots, \lambda) = 0. \tag{A'}
\end{equation}
\begin{equation}\label{B'}
\pd{g}{z_j}(\lambda,\ldots, \lambda) = 0. \tag{B'}
\end{equation}
\begin{equation}\label{C'}
g(z_1+t,\ldots,z_n+t) = g(z_1,\ldots, z_n),\mbox{ for all } t\in \R. \tag{C'}
\end{equation}
Condition $\ref{C'}$ implies that 
\begin{equation}\label{C''}
g(z_1+c, \ldots, z_n+c) = g(z_1,\ldots, z_n),
\end{equation}
for all complex $c$ with $\im(c) > - \min_j \im(z_j) $.
Indeed, fixing $z_1,\ldots, z_n\in \H$, the holomorphic function $c \mapsto g(z_1+c, \ldots, z_n+c) - g(z_1,\ldots, z_n)$ vanishes
on the real axis and, thus, on the whole domain. 

Now, write $g(z_1,\ldots, z_n) = h(a, d_1, \ldots, d_{n-1})$,
where $$a = \frac{1}{n}\sum_{j=1}^n z_j\mbox{ and } d_j = z_{j+1}-z_j \mbox{ for } j=1,\ldots, n-1,$$
and $h$ is holomorphic on the image $\Omega$ of $\H^n$ under the coordinate change.

For each $a\in \H$, let $$\Omega(a)=\{(d_1,\ldots, d_n) \in \C^{n-1}~|~(a,d_1,\ldots, d_n) \in \Omega\}.$$
Define $h^a: \Omega(a) \ra \C$ by
$$h^a(d_1,\ldots, d_{n-1}) = h(a, d_1,\ldots, d_{n-1}).$$

For each $a\in \H$, $\Omega(a)$ is a convex open set containing the origin. Moreover, $\Omega(ta) = t\Omega(a)$ for $t>0$. It follows that $\Omega(i t_1)\subset \Omega(i t_2)$ for $0<t_1< t_2$,
and that $\bigcup_{t>0} \Omega(i t) = \C^{n-1}$.

Now, \eqref{C''} implies $h^{it_1}(d_1,\ldots, d_{n-1}) = h^{it_2}(d_1,\ldots, d_{n-1})$ whenever $(d_1,\ldots, d_{n-1}) \in \Omega(i t_1)$
and $t_1<t_2$. Since $\bigcup_{t>0} \Omega(i t) = \C^{n-1}$, there is a holomorphic $H: \C^{n-1} \ra \C$ so that
$h^{it} = H|_{\Omega(i t)}.$ Again by \eqref{C''}, $h^{x+iy} = h^{iy}$ for all $x+iy\in \H$. So
$$h^{a} = H|_{\Omega(a)},~\forall a\in \H.$$ 
It thus suffices to show that $H$ is identically 0.

Recall that
$$f = \frac{1}{n}\sum_{j=1}^{n} z_j + g(z_1,\ldots, z_n) = a + h(a,d_1,\ldots, d_{n-1}).$$
Since $f$ maps into $\H$, $h^a$ maps $\Omega(a)$ into the strip
$\{z~|~\im(z) > -\im(a) \}.$
Thus, $H$ maps each $\Omega(i t)$ to $\{\im(z) > -t\}$.

Recall that $\Omega(i)$ is open and contains 0. So $\Omega(i)$ contains an open Euclidean ball $B_r(0)$ centered at the origin. Then $B_{rt}(0) \subset \Omega(i t)$,
so $H(B_{rt}) \subset \{\im(z)>-t\}$, for all $t>0$. Thus, if $\lp {\sum |d_j|^2}\rp^{1/2}=rt$, then $$\im \lbrac H(d_1,\ldots, d_{n-1})\rbrac \geq -t.$$
In other words, we have
\begin{equation}\label{One}
\im\lbrac H(d_1,\ldots, d_{n-1}) \rbrac \geq -\frac{\lp {\sum |d_j|^2}\rp^{1/2}}{r},
\end{equation} 
for all $(d_1,\ldots d_{n-1})\in \C^{n-1}$.

Condition \eqref{A'} implies $h(a,0,\ldots, 0)=0$, so
\begin{equation}\label{Two}
H(0,\ldots,0) = 0.
\end{equation}

Finally, condition \eqref{B'} and the chain rule imply that the derivatives
$\pd{h}{d_j}(a,0,\ldots,0)$ are 0, so that
\begin{equation}\label{Three}
\pd{H}{d_j}(0,\ldots,0) = 0~ \forall j.
\end{equation}

We reduce to Lemma \ref{Harmonic}. Fix arbitrary $(d_1,\ldots,d_{n-1})$ with
$\sum_j |d_j|^2 = 1$.
By $\eqref{One},\eqref{Two}$, and $\eqref{Three}$ the harmonic function $$\phi(z) = \im\lbrac H(d_1z, \ldots, d_{n-1} z ) \rbrac$$
satisfies the conditions of the lemma, with $C = \frac{1}{r}$.
We conclude that $\im(H)$, and thus $H$, are identically 0.
\qed 

\section{Appendix: Polarization}\label{Seven}

Markovic's proof in \cite{Ma} of Proposition \ref{Ma} uses the classical polarization principle.
The proof generalizes almost verbatim to a proof of the corresponding result for maps $\H^n \ra \H$ (Corollary \ref{Cor}), but the polarization principle must be 
replaced by the following fact: 
\begin{Proposition}\label{Q}
Let $V$ be the real vector subspace of $\C^n$ consisting
of points the form $(r + t_1 i, \ldots, r + t_n i)$ with $r$ and $t_1,\ldots, t_n$ real and $\sum_{j=1}^n t_j=0$. 
Let $U\subset \C^n$ be a domain such that $U\cap V$ is nonempty.
If $h:U\ra \C$ is holomorphic and vanishes on $U\cap V$, then $h$ is identically 0 on $U$. 
\end{Proposition}
 
\noindent (The polarization principle is the $n=2$ case of the above result.) We will prove Proposition \ref{Q} as a corollary of the following well-known proposition.
\begin{Proposition}\label{P}
Let $U\subset \C^n$ be a domain, and let $M\subset U$ be a nonempty smooth submanifold. Suppose for each $p\in M$ that
 $T_pM$ and $i\lp T_pM \rp$ together span $\C^n$. Let $h:U\ra \C$ be a holomorphic function which vanishes on $M$. 
 Then $h$ is identically 0 on $U$.
\end{Proposition}

\noindent{\em Proof:}
 Let $p\in M$, and consider the differential $dh_p: \C^n \ra \C$. Since $f$ vanishes on $M$, $dh_p$ vanishes on $T_pM$.
Since $dh_p$ is complex-linear, it vanishes also on $i\lp T_pM\rp$. But since $T_pM + i\lp T_pM \rp=\C^n$,
$dh_p=0$. Since $p$ was arbitrary, we conclude the first partial derivatives $\pd{h}{z_j}$ vanish on $M$.
Applying the same argument to $\pd{h}{z_j}$, we find that the second partials $\frac{\partial^2h}{\partial z_k\partial z_j}$ also vanish on $M$.
Continuing inductively, we find that all higher derivatives vanish on $M$. Since $h$ is analytic, it follows that $h$ is identically 0 on $U$. 
\qed

\vspace{2mm}

\noindent{\em Proof of Proposition \ref{Q}:}
If $p \in U\cap V$, $T_pV$ identifies naturally with $V$. The vector space $V$ has (real) dimension $n$, and $V \cap iV = \{0\}$,
so $\C^n = V \oplus iV$. So Proposition \ref{P} applies, with $M = U\cap V$. 

\qed   

\section{Acknowledgments}
I would like to thank Peter Burton, Oleg Ivrii, Gregory Knese, John M\textsuperscript{c}Carthy, and Vladimir Markovic for helpful discussion.

\end{document}